\documentclass[a4paper,american,reqno]{amsart}

\usepackage{babel}
\usepackage[utf8]{inputenc}
\usepackage[T1]{fontenc}
\usepackage{xspace}
\usepackage{todonotes}
\usepackage{csquotes}
\usepackage{amssymb}
\usepackage{microtype}
\usepackage{pgfplots}
\usepackage{tikz}
\usepackage[style=authoryear-comp,
            maxbibnames=100,
            maxcitenames=2,
            dashed=false,
            uniquename=false,
            giveninits=true,
            isbn=false,
            sorting=nty,
            backend=biber]{biblatex}
\usepackage[colorlinks,
            citecolor=blue,
            urlcolor=blue,
            linkcolor=blue]{hyperref} 

\usetikzlibrary{plotmarks}
\usepgfplotslibrary{fillbetween}
\pgfplotsset{compat=1.16}

\makeatletter
\patchcmd{\@settitle}{\uppercasenonmath\@title}{\scshape\large}{}{}
\patchcmd{\@setauthors}{\MakeUppercase}{\scshape\normalsize}{}{}
\makeatother

\makeatletter
\@namedef{subjclassname@2020}{%
  \textup{2020} Mathematics Subject Classification}
\makeatother

\vfuzz \hfuzz
\newcommand{\st}{\text{s.t.}}
\newcommand{\field}{\mathbb}

\newcommand{\reals}{\field{R}}
\newcommand{\integers}{\field{Z}}

\newcommand{\R}{\reals}
\newcommand{\Z}{\integers}

\DeclareMathOperator*{\argmin}{arg\,min}
\DeclareMathOperator*{\argmax}{arg\,max}

\newcommand{\set}[1]{\{#1\}}
\newcommand{\Set}[1]{\left\{#1\right\}}
\newcommand{\defset}[3][\defsep]{\set{#2#1#3}}
\newcommand{\Defset}[3][\defsep]{\Set{#2#1#3}}

\newcommand{\abbr}[1][abbrev]{#1.\xspace}

\newcommand{\eg}{\abbr[e.g]}

\newcommand{\ie}{\abbr[i.e]}

\newcommand{\wrt}{\abbr[w.r.t]}

\newcommand{\define}{\mathrel{{\mathop:}{=}}}

\definecolor{my-royal-blue}{HTML}{2250d9}
\definecolor{my-complementary-orange}{HTML}{d94f22}

\definecolor{ukrainian-blue}{RGB}{0,87,183}
\definecolor{ukrainian-yellow}{RGB}{255,215,0}

\begin{document}

\title{A Brief Introduction to Robust Bilevel Optimization}

\author[Y. Beck]{Yasmine Beck}
\author[I. Ljubi\'c]{Ivana Ljubi\'c}
\author[M. Schmidt]{Martin Schmidt}

\address[Y. Beck, M. Schmidt]{%
  Trier University,
  Department of Mathematics,
  Universitätsring 15,
  54296 Trier,
  Germany}%
\email{yasmine.beck@uni-trier.de}
\email{martin.schmidt@uni-trier.de}

\address[I. Ljubi\'c]{%
  ESSEC Business School of Paris,
  Cergy-Pontoise,
  France}%
\email{ljubic@essec.edu}

\date{\today}

\begin{abstract}
  Bilevel optimization is a powerful tool for modeling hierarchical
decision making processes.
However, the resulting problems are challenging to solve---both in
theory and practice.
Fortunately, there have been significant algorithmic advances in the
field so that we can solve much larger and also more complicated
problems today compared to what was possible to solve two decades
ago.
This results in more and more challenging bilevel problems that
researchers try to solve today.
In this article, we give a brief introduction to one of these more
challenging classes of bilevel problems:
bilevel optimization under uncertainty using robust optimization techniques.
To this end, we brief\/ly state different versions of uncertain bilevel
problems that result from different levels of cooperation of the
follower as well as on when the uncertainty is revealed.
We highlight these concepts using an academic example and discuss
recent results from the literature concerning complexity as well as
solution approaches.
Finally, we discuss that the sources of uncertainty in bilevel
optimization are much richer than in single-level optimization and, to
this end, introduce the concept of decision uncertainty.


\end{abstract}


\maketitle

\section{Introduction}
\label{sec:introduction}

Bilevel optimization has its roots in economics and dates back to the
seminal works by \textcite{Stackelberg:1934,Stackelberg:1952}.
It has been introduced in the field of mathematical optimization much
later in the publications by \textcite{Bracken-Mcgill:1973} as well as
\textcite{Candler-Norton:1977}.
We use bilevel optimization to model hierarchical decision making processes,
typically with two players, which we refer to as the leader and the follower.
Despite its intrinsic hardness \parencite{jeroslow:1985,Hansen-et-al:1992},
several innovative works pushed the boundaries of computational bilevel
optimization so that we can tackle some relevant practical
applications today; see, \eg, \textcite{Kleinert_et_al:2021} for a
recent survey on computational bilevel optimization as well as the
annotated bibliography by \textcite{Dempe:2020}.

The main goal of this article is to give a brief introduction to some
basic concepts of bilevel optimization problems under uncertainty.
The field is still in its infancy but, nevertheless, due to its
relevance in many practical applications, it is developing very fast.
In classic, \ie, single-level, optimization, there are two major
approaches to address uncertainty:
stochastic optimization \parencite{Birge_Louveaux:2011,
  Kall_Wallace:1994} and robust optimization \parencite{Soyster:1973,
  Ben-Tal_Nemirovski:1998,Ben-Tal_et_al:2009,Bertsimas_et_al:2011}.
The same two paths have been followed as well in bilevel optimization
starting from the 1990s on.
However, the sources of uncertainty are much richer in bilevel
optimization compared to single-level optimization.
To make this more concrete, let us consider the linear optimization
problem $\min_x \defset{c^\top x}{Ax \geq b}$.
It can ``only'' be subject to uncertainty due to uncertainties in the
problem's data $c$, $A$, and $b$.
Throughout this article, we will refer to this setting as \emph{data
  uncertainty}.
Moreover, a bilevel optimization problem may also be subject to an
additional source of uncertainty, which is due to its nature that
it combines two different decision makers in one model.
Hence, there can be further uncertainty involved either if the leader
is not sure about the reaction of the follower or if the follower is
not certain about the observed leader's decision.
We will denote this additional type of uncertainty as \emph{decision
  uncertainty}.
Obviously, decision uncertainty does not play any role in single-level
optimization since only one decision maker is involved.

In this introductory article, we will solely focus on data uncertainty
that is tackled using concepts from robust optimization.
For more details regarding stochastic bilevel optimization, decision
uncertainty, etc.\ we refer to our recent
survey \parencite{Beck_et_al:2022c}.


\section{Problem Statement}
\label{sec:problem-statement}

We start by considering the deterministic bilevel problem (we explain
the quotation marks below)
\begin{subequations}
  \label{eq:general-UL}
  \begin{align}
    ``\min_{x \in X}" \quad
    & F(x,y) \label{eq:general-UL-obj} \\
    \st \quad
    & G(x,y) \geq 0, \\
    & y \in S(x),
  \end{align}
\end{subequations}
where~$S(x)$ denotes the set of optimal solutions of
the~$x$-parameterized problem
\begin{subequations} \label{eq:general-LL}
  \begin{align}
    \min_{y \in Y} \quad & f(x,y) \label{eq:general-LL-obj} \\
    \st \quad & g(x,y) \geq 0.
  \end{align}
\end{subequations}
Problem~\eqref{eq:general-UL} is referred to as the upper-level
(or the leader's) problem and Problem~\eqref{eq:general-LL} is
the so-called lower-level (or the follower's) problem.
Moreover, we refer to~$x \in X$ and~$y \in Y$ as the leader's and the
follower's variables, respectively.
The sets~\mbox{$X \subseteq \R^{n_x}$} and~$Y \subseteq \R^{n_y}$ can be
used to include possible integrality constraints.
The objective functions are given by~$F, f : \R^{n_x} \times \R^{n_y} \to \R$
and the constraint functions by~\mbox{$G : \R^{n_x} \times \R^{n_y} \to \R^{m}$}
as well as~$g : \R^{n_x} \times \R^{n_y} \to \R^\ell$.
In the case that the lower-level problem does not have a unique solution,
the bilevel problem~\eqref{eq:general-UL} and~\eqref{eq:general-LL} is
ill-posed. This ambiguity is expressed by the quotation marks
in~\eqref{eq:general-UL-obj}. To overcome this issue, it is common to pursue
either an optimistic or a pessimistic approach to bilevel optimization;
see, \eg, \textcite{Dempe:2002}.
In the optimistic setting, the leader chooses the follower's response
among multiple optimal solutions of the lower-level problem such
that it favors the leader's objective function value.
Hence, the leader also minimizes her\footnote{Throughout this article,
  we use ``her'' for the leader and ``his'' for the follower.}
objective in the~$y$ variables, \ie, we consider the problem
\begin{equation} \label{eq:general-optimistic-bilevel}
  \min_{x \in \bar{X}} \ \min_{y \in S(x)} \quad F(x,y)
\end{equation}
with~$\bar{X} \define \defset{x \in X}{G(x) \geq 0}$
and~$G : \R^{n_x} \to \R^m$.
Here and in what follows, we focus on the setting without coupling constraints,
\ie, without upper-level constraints that depend on the
variables~$y$.
In the pessimistic setting, the leader anticipates that, among multiple
optimal solutions of the follower, the worst possible response \wrt \ the
upper-level objective function will be chosen by the follower.
Thus, one studies the problem
\begin{equation*}
  \min_{x \in \bar{X}} \ \max_{y \in S(x)} \quad F(x,y).
\end{equation*}

In this article, we focus on bilevel problems of the above form which
are additionally affected by data uncertainty.

\subsection{Data Uncertainty}
\label{sec:data-uncertainty}

Data uncertainty arises when some of the players only have access
to inaccurate or incomplete data.
In robust optimization, it is assumed that these uncertainties take
values in a given, and usually compact, uncertainty set~$\mathcal{U}$.
The uncertainty sets are typically modeled using boxes,
polyhedra, ellipsoids, or cones; see, \eg,
\textcite{Bertsimas_et_al:2011,Ben-Tal_Nemirovski:1998,Ben-Tal_et_al:2004,
  Ben-Tal_et_al:2009,Soyster:1973}.
In the context of single-level robust optimization, there are two
possibilities to hedge against data uncertainty.

First, assuming that the coefficients of the
objective function are uncertain, one searches for a solution that is
optimal for the worst-case realization of the uncertain parameters.
The problem can be modeled as
\begin{equation} \label{eq:single-level-robust-obj}
  \min_{x \in \bar{X}} \ \max_{u  \in \mathcal{U}} \quad F(x,u),
\end{equation}
where the objective function $F : \R^{n_x} \times \R^{n_u} \to \R$
is continuous and the sets~$\mathcal{U} \subseteq
  \R^{n_u}$ and~$\bar{X}$ are defined as above.

Second, in the case that the uncertainty affects the
coefficients of the constraints, one is interested in a solution that
is feasible for all possible realizations of the uncertainty.
This problem can be stated as
\begin{equation} \label{eq:single-level-robust-constr}
  \min_{x \in X} \quad F(x)
  \quad \st \quad
  G(x,u) \geq 0 \quad \text{for all } u \in \mathcal{U},
\end{equation}
where both the objective function~$F : \R^{n_x} \to \R$
and the constraint function~\mbox{$G : \R^{n_x} \times \R^{n_u} \to \R^{m}$}
are continuous.
Problem~\eqref{eq:single-level-robust-constr} can be reformulated as
\begin{equation} \label{eq:single-level-robust-constr-min}
  \min_{x \in X} \quad F(x) \quad
  \st \quad \min \Defset{ G(x,u) }{ u \in \mathcal{U} } \geq 0.
\end{equation}
In particular, Problem~\eqref{eq:single-level-robust-obj} can be restated
as an instance of Problem~\eqref{eq:single-level-robust-constr-min} using an
epigraph reformulation, \ie,
\begin{equation*}
  \min_{x \in \bar{X}, t \in \R} \quad t
  \quad \st \quad t \geq \max \Defset{ F(x,u) }{ {u \in \mathcal{U}} }.
\end{equation*}

Note that for the two settings discussed so far, a single decision
maker has to take a here-and-now decision before the uncertainty is
revealed.
In bilevel optimization, however, there are two different timings that
are possible---one in which the uncertainty realizes after and one in
which the uncertainty realizes before the follower takes his decision.

\subsubsection{Here-and-Now Follower}
\label{sec:here-now-follower}

In this case, both the leader
and the follower have to make their decisions before the
uncertainty is revealed, \ie, one considers the timing
\begin{equation} \label{eq:here-and-now-follower-timing}
  \text{leader } x
  \quad\curvearrowright\quad
  \text{follower } y = y(x)
  \quad\curvearrowright\quad
  \text{uncertainty } u.
\end{equation}
This means that the leader anticipates an optimal response of the
follower who hedges against data uncertainty.
Hence, the lower-level problem is an~$x$-parameterized problem in which
we can embed any of the concepts known for single-level optimization
under uncertainty.
For instance, if only the lower-level objective function is uncertain
and the follower is assumed to behave in an optimistic way,
we are solving Problem~\eqref{eq:general-optimistic-bilevel}
with
\begin{equation*}
  S(x) \define \argmin_{y' \in Y} \Defset{
    \max_{u \in \mathcal{U}} \ f(x,u,y')
  }{
    g(x,y') \geq 0}.
\end{equation*}

\subsubsection{Wait-and-See Follower}
\label{sec:wait-see-follower}

In this setting, the leader first
takes a here-and-now decision, \ie,
without knowing the realization of uncertainty.
Then, the uncertainty is revealed and,
finally, the follower decides in a wait-and-see fashion,
taking the leader's decision as well as the realization of the uncertainty
into account. Hence, one considers the timing
\begin{equation} \label{eq:wait-and-see-follower-timing}
  \text{leader } x
  \quad\curvearrowright\quad
  \text{uncertainty } u
  \quad\curvearrowright\quad
  \text{follower } y = y(x,u).
\end{equation}
This means that the leader does not have full knowledge about the
lower-level problem.
Thus, she wants to hedge against the worst-case reaction of the
follower.
Here, ``worst-case'' may not only imply the robustness
of the leader \wrt \ lower-level data uncertainty but also
her conservatism regarding the cooperation of the follower.
For instance, to protect against the worst-case realization of
the uncertainties \wrt \ the leader's objective function, we consider
the problem
\begin{equation}
  \label{eq:general-UL-u-OF}
  ``\min_{x \in \bar{X}} \ \max_{u \in \mathcal{U}}" \quad F(x,y)
  \quad \st \quad y \in S(x,u),
\end{equation}
where~$S(x,u)$ is the set of optimal solutions of
the~$(x,u)$-parameterized problem
\label{eq:general-LL-u}
\begin{equation*}
  \min_{y \in Y} \quad f(x,u,y)
  \quad \st \quad  g(x,u,y) \geq 0.
\end{equation*}
The quotation marks in~\eqref{eq:general-UL-u-OF} express the
ill-posedness of the bilevel problem in the case that the set~$S(x,u)$
is not a singleton.
Hence, one also needs to distinguish between the optimistic and
the pessimistic case in the robust setting.
Indeed, both situations can be motivated by practical applications.
For instance, the pessimistic robust bilevel problem appears when the
leader wants to hedge against the worst-case both w.r.t.\
lower-level data uncertainty as well as w.r.t.\ the potentially
unknown level of cooperation of the follower.
On the other hand, there may also be situations in which the follower
still hedges against his uncertainties in a robust way but, in the case
of ambiguous optimal solutions, acts in an optimistic way.
This might be the case in energy markets with sufficiently regulated
agents, where a strong level of regulation might lead to an optimistic
robust bilevel problem.


\section{An Academic Example}
\label{sec:example}

\begin{figure}
  \begin{center}
    \newcommand*\obj[1]{4.75 - #1}
\begin{tikzpicture}[scale=0.7]
	\begin{axis}[
    axis lines = left,
    ymin = 0,
    ymax = 3.5,
    xmin = 0,
    xmax = 4.5,
    xtick = {0,1,2,3,4},
    xticklabels = {,1,2,3,4},
    ytick = {0,1,2,3},
	  yticklabels = {,1,2,3},
	  xlabel = {$x$},
	  ylabel = {$y$},
	  every axis x label/.style={
  	  at={(ticklabel* cs:1)},
    	  anchor=west,
	  },
	  every axis y label/.style={at={(ticklabel* cs:1)},anchor=south,},
    clip = false
    ]
	\addplot [name path=xaxis,
	    domain=1:3.5,
	    samples=2,
	    color=black,
	    very thick,
	    ]{0};
	\addplot [name path=parallel,
	    domain=1.5:3,
	    samples=2,
	    color=my-complementary-orange,
	    very thick,
	    ]{2.5};
    \addplot [
        domain=0.5:1,
        samples=2,
        color=my-royal-blue,
	    dashed,
	    very thick,
        ]{3-3*x};
    \addplot [name path=b,
        domain=3.5:4,
        samples=2,
        color=black,
	    very thick,
        ]{-7+2*x};
    \addplot [name path=a,
        domain=0.5:1.5,
        samples=2,
        color=my-royal-blue,
	    dashed,
	    very thick,
        ]{1+x};
    \addplot [
        domain=3:4,
        samples=2,
        color=my-complementary-orange,
	    very thick,
        ]{7-1.5*x};
    \addplot[ultra thick,black,mark=*,mark size = 2] coordinates {(1.5,2.5)};
	
	\addplot[color=gray!50, opacity=0.3] fill between [of = a and xaxis, soft clip = {domain=0.5:4}];
    \addplot[color=gray!50, opacity=0.3] fill between [of = b and parallel, soft clip = {domain=1.5:4}];
	
	\addplot[thick, dashed, color=my-royal-blue, domain=1.146:1.854]{\obj{\x}};
	\coordinate (p1) at (axis cs:1.5,3.25);
	\coordinate (p2) at (axis cs:1.323,3.073);
	\draw[thick, color=my-royal-blue, ->, style = {-to}] (p1) -- (p2);
	\draw[] (p1) node[anchor = south west, align = center, color=my-royal-blue] {$F$};

	\addplot[thick, domain=3.5:4.5]{2.5};
	\coordinate (p5) at (axis cs:4,2.5);
	\coordinate (p6) at (axis cs:4,2.75);
	\draw[thick, ->, style = {-to}] (p5) -- (p6);
        \draw[] (p5) node[anchor = north, align = center] {$f$};
	\end{axis}
\end{tikzpicture}
    \quad
    \begin{tikzpicture}[scale=0.7]
	\begin{axis}[
    axis lines = left,
    ymin = 0,
    ymax = 3.5,
    xmin = 0,
    xmax = 4.5,
    xtick = {0,1,2,3,4},
    xticklabels = {,1,2,3,4},
    ytick = {0,1,2,3},
	  yticklabels = {,1,2,3},
	  xlabel = {$x$},
	  ylabel = {$y$},
	  every axis x label/.style={
  	  at={(ticklabel* cs:1)},
    	  anchor=west,
	  },
	  every axis y label/.style={at={(ticklabel* cs:1)},anchor=south,},
    clip = false
    ]
	\addplot [name path=xaxis,
	       domain=1:3.5,
	       samples=2,
	       color=my-complementary-orange,
	       very thick,
	       ]{0};
	\addplot [name path=parallel,
	       domain=1.5:3,
	       samples=2,
	       color=black,
	       very thick,
	       ]{2.5};
    \addplot [
          domain=0.5:1,
          samples=2,
          color=my-royal-blue,
	       dashed,
	       very thick,
          ]{3-3*x};
    \addplot [name path=b,
          domain=3.5:4,
          samples=2,
          color=my-complementary-orange,
	       very thick,
          ]{-7+2*x};
    \addplot [name path=a,
          domain=0.5:1.5,
          samples=2,
          color=my-royal-blue,
	       dashed,
	       very thick,
          ]{1+x};
    \addplot [
          domain=3:4,
          samples=2,
          color=black,
	       very thick,
          ]{7-1.5*x};
    \addplot[ultra thick,black,mark=*,mark size = 2] coordinates {(1,0)};
	
	\addplot[color=gray!50, opacity=0.3] fill between [of = a and xaxis, soft clip = {domain=0.5:4}];
    \addplot[color=gray!50, opacity=0.3] fill between [of = b and parallel, soft clip = {domain=1.5:4}];
	
	\addplot[
        domain=1.146:1.854,
        samples=2,
        thick,
        dashed,
        color=my-royal-blue,
        ]{4.75 - x};
	\coordinate (p1) at (axis cs:1.5,3.25);
	\coordinate (p2) at (axis cs:1.323,3.073);
	\draw[thick, color=my-royal-blue, ->, style = {-to}] (p1) -- (p2);
	\draw[] (p1) node[anchor = south west, align = center, color=my-royal-blue] {$F$};

	\addplot[thick, domain=3.5:4.5]{2.5};
	\coordinate (p5) at (axis cs:4,2.5);
	\coordinate (p6) at (axis cs:4,2.25);
	\draw[thick, ->, style = {-to}] (p5) -- (p6);
        \draw[] (p5) node[anchor = south, align = center] {$\tilde{f}$};
	\end{axis}
\end{tikzpicture}
    \caption{Both figures show the upper-level constraints (dashed blue
      lines), the lower-level constraints (solid black and orange
      lines), the shared constraint set (gray area), and the bilevel
      feasible set (solid orange lines) of the bilevel
      problem~\eqref{eq:example-UL} and~\eqref{eq:example-LL}.
      The deterministic variant of the problem is depicted on the left
      and the variant with a here-and-now follower is given on the
      right.}
    \label{fig:example}
  \end{center}
\end{figure}
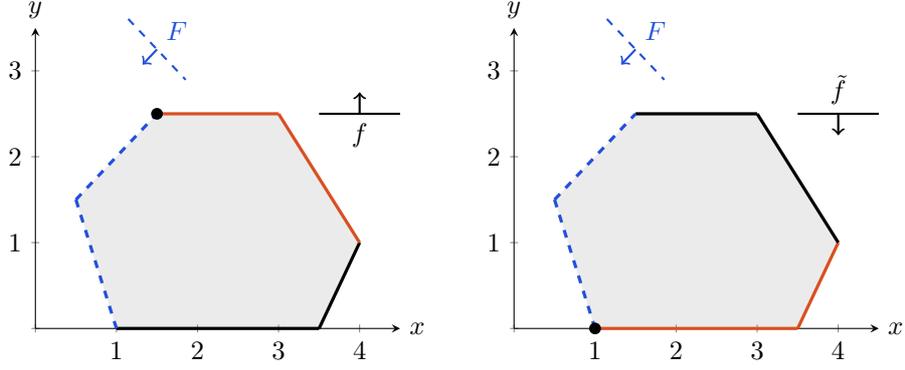

Let us consider the linear bilevel problem taken from
\textcite{Beck_et_al:2022c} that is given by
\begin{subequations}
  \label{eq:example-UL}
  \begin{align}
    ``\min_{x \in \R}" \quad
    & F(x,y) = x + y \\
    \st \quad
    \ & x - y \geq -1,
        \label{eq:example-UL-1} \\
    & 3x + y \geq 3,
      \label{eq:example-UL-2}\\
    & y \in S(x),
  \end{align}
\end{subequations}
where~$S(x)$ denotes the set of optimal solutions of the~$x$-parameterized
lower level
\begin{subequations} \label{eq:example-LL}
  \begin{align}
    \min_{y \in \R} \quad
    & f(x, y) = -0.1 y \\
    \st \quad
    & -2x + y \geq -7, \label{eq:example-LL-1}\\
    & -3x - 2y \geq -14, \label{eq:example-LL-2}\\
    & 0 \leq y \leq 2.5. \label{eq:example-LL-3}
  \end{align}
\end{subequations}
The problem is depicted in Figure~\ref{fig:example}~(left).
The upper- and lower-level constraints are represented with dashed and
solid lines, respectively.
The optimal solution~$(x^*,y^*)=(1.5,2.5)$ is the same for both the optimistic
and the pessimistic setting and it is illustrated by the thick
dot.
Suppose now that the lower-level objective function is uncertain.
To this end, we consider~$\tilde{f}(x,u,y) = (-0.1 + u)y$ and
assume that~$u$ only takes values in the uncertainty
set~$\mathcal{U} = \defset{u \in \R}{|u| \leq 0.5}$.
In what follows, we distinguish between a follower taking a
here-and-now or a wait-and-see decision to illustrate how the
considered timing may affect the solution of the problem.

\subsection{Here-and-Now Follower}

We first consider the timing
in~\eqref{eq:here-and-now-follower-timing}.
The robustified lower-level problem is thus given by
\begin{equation*}
  \min_{y \in \R} \ \max_{u \in \mathcal{U}}
  \quad \tilde{f}(x,u,y) = (-0.1 + u)y
  \quad \st \quad
  \text{\eqref{eq:example-LL-1}--\eqref{eq:example-LL-3}}.
\end{equation*}
Using classic techniques from robust optimization, we obtain a
modified gradient of the lower-level objective function, which is
shown in Figure~\ref{fig:example}~(right).
The optimal solution~\mbox{$(x^*,y^*)=(1,0)$} of this problem is
represented by the thick dot.
In particular, there is a unique lower-level response for every
feasible~$x$, which is why we do not need to distinguish between
the optimistic and the pessimistic case.

\subsection{Wait-and-See Follower}

We now consider the timing
in~\eqref{eq:wait-and-see-follower-timing}, \ie, the overall
robustified bilevel problem reads
\begin{equation*}
  ``\min_{x \in \R} \ \max_{u \in \mathcal{U}}" \quad F(x,y)
  \quad \st \quad
  \text{\eqref{eq:example-UL-1}--\eqref{eq:example-UL-2}},
  \ y \in S(x,u),
\end{equation*}
where~$S(x,u)$ is the set of optimal solutions of the
$(x,u)$-parameterized lower level
\begin{equation*}
  \min_{y \in \R} \quad \tilde{f}(x,u,y) = (-0.1 + u)y
  \quad \st \quad
  \text{\eqref{eq:example-LL-1}--\eqref{eq:example-LL-3}}.
\end{equation*}
To solve this problem, we need to distinguish the following three
cases.
\begin{enumerate}
\item $-0.5 \leq u < 0.1$:
  This case corresponds to the setting that is
  depicted in Figure~\ref{fig:example}~(left).
  The optimal follower's reaction is thus given by
  \begin{equation}
    \label{eq:example-worst-follower-reaction}
    y(x, u) =
    \begin{cases}
      2.5, & x \leq 3,\\
      -1.5x + 7, & 3 \leq x \leq 4.
    \end{cases}
  \end{equation}
  Note, however, that the bilevel problem is infeasible for~$x <
  1.5$.
  In particular, this means that the robust optimal leader's
  decision~$x^*=1$ for the case with a here-and-now follower
  is no longer bilevel feasible if the follower decides in a
  wait-and-see fashion.
\item $u = 0.1$:
  Any feasible decision of the follower, \ie, any~$y \in \R$ that
  satisfies~\eqref{eq:example-LL-1}--\eqref{eq:example-LL-3}, is optimal
  for the~$x$-parameterized lower level.
  Hence, the distinction between an optimistic and
  a pessimistic follower is necessary.
  In the optimistic setting, the follower would react with
  \begin{equation}
    \label{eq:example-best-follower-reaction}
    y(x, u) =
    \begin{cases}
      0, & x \leq 3.5,\\
      2x - 7, & 3.5 \leq x \leq 4.
    \end{cases}
  \end{equation}
  This corresponds to the setting that is depicted
  in Figure~\ref{fig:example}~(right).
  A pessimistic follower, however, would
  select~\eqref{eq:example-worst-follower-reaction}.
  Note that the bilevel problem with an optimistic
  follower turns out to be infeasible for~$x < 1$ and, again, the
  problem is infeasible for~$x < 1.5$ if a pessimistic follower is
  considered.
\item $0.1 < u \leq 0.5$:
  The optimal follower's reaction is given
  by~\eqref{eq:example-best-follower-reaction}.
  Again, the overall bilevel problem turns out to be infeasible
  for~$x < 1$.
\end{enumerate}

To determine an optimal solution of the bilevel
problem~\eqref{eq:example-UL} and~\eqref{eq:example-LL} with a
wait-and-see follower, we thus consider the worst-case realization
of each of the previous three cases \wrt\ the leader's
decision~$x$.
Hence, we need to solve
\begin{equation} \label{eq:example-wait-and-see}
  \min_x \quad \hat{F}(x) \quad \st \quad 1.5 \leq x \leq 4
\end{equation}
with the piecewise-linear function
\begin{equation*}
  \hat{F}(x) =
  \begin{cases}
    x + 2.5,
    & 1.5 \leq x \leq 3,\\
    -0.5x + 7,
    & 3 \leq x \leq 4.
  \end{cases}
\end{equation*}
In particular, the solution~$x^* = 1.5$ of
Problem~\eqref{eq:example-wait-and-see}
is an optimal decision of the leader in both the optimistic and the
pessimistic setting. After observing the realization
of the uncertainty, the corresponding response of the follower is then
given by
\begin{equation*}
  y^*_{\text{o}}(x^*,u) =
  \begin{cases}
    2.5, & -0.5 \leq u < 0.1,\\
    0, & 0.1 \leq u \leq 0.5
  \end{cases}
\end{equation*}
in the optimistic setting, whereas, for the pessimistic case, we have
\begin{equation*}
  y^*_{\text{p}}(x^*,u) =
  \begin{cases}
    2.5, & -0.5 \leq u \leq 0.1,\\
    0, & 0.1 < u \leq 0.5.
  \end{cases}
\end{equation*}
Note that, if~$u \in [-0.5, 0.1)$ realizes, at the point $x^* = 1.5$, the deterministic
solution $(x^*,y(x^*))$ and the robust bilevel solutions $(x^*,y(x^*,u))$ coincide. However,
the optimal follower's response $y(x^*,u)$ in the robust setting may change
significantly for~$u \geq 0.1$.


\section{Selected Results from the Literature}
\label{sec:recent-literature-results}

The field of robust bilevel optimization is still in its infancy.
For a detailed discussion of existing modeling and solution approaches,
we refer to our recent survey \parencite{Beck_et_al:2022c}.
In deterministic bilevel optimization, a standard solution approach is
to reformulate the problem as a classic, \ie, single-level, problem.
This can be done, \eg, by replacing the lower level with its
Karush--Kuhn--Tucker (KKT) conditions \parencite{Fortuny-Amat-Mccarl:1981}.
The same holds true for robust bilevel problems whenever the
robust counterpart of the lower-level problem can be reformulated as a
deterministic problem for which the KKT conditions are necessary and
sufficient.
However, these reformulation techniques cannot be applied
anymore if discrete variables are introduced in the lower level.
Due to their intrinsic hardness, approaches for discrete robust
bilevel problems have not been investigated a lot up to now.
In single-level optimization, the knapsack problem is one of the most
thoroughly studied discrete optimization problem due to its relevance both
in theory and practice; see, \eg, \textcite{Pisinger_Toth:1998}.
Bilevel knapsack problems naturally extend their single-level counterparts
such as to capture hierarchical and, in particular, competitive settings
\parencite{DeNegre:2011,Caprara_et_al:2013,Fischetti_et_al:2018b,
  Fischetti_et_al:2019,Della-Croce_Scatamacchia:2020}.
Moreover, the bilevel knapsack interdiction problem is commonly
used as a benchmark for testing bilevel optimization solvers; see, \eg,
\textcite{DeNegre_Ralphs:2009,Tang_et_al:2016}.
It is thus not surprising that bilevel knapsack problems are also among
the first discrete bilevel problems studied under uncertainty---both in terms of
complexity questions and solution approaches.
The remainder of this section is thus dedicated to a brief overview of
recent results from the literature for robust bilevel knapsack problems.

\subsection{Complexity Results for Robust Continuous Bilevel Knapsack Problems with a Wait-and-See Follower}

We start by considering the robust continuous bilevel knapsack problem
with an uncertain lower-level objective, \ie, we consider the problem
\begin{subequations} \label{eq:continuous-robust-knapsack}
  \begin{align}
    \max_{x \in [x^-,x^+]} \ \min_{c \in \mathcal{U},y \in \R^n} \quad & d^\top y
    \\
    \st \quad  & y \in \argmax_{y'} \Defset{c^\top y'
                 }{
                 a^\top y' \leq x,\, 0 \leq y' \leq 1}
  \end{align}
\end{subequations}
with~$x^-,\, x^+ \in \R$, $x^- \leq x^+$, $a,\, c,\, d \in \R^n_{\geq 0}$,
and an uncertainty set~$\mathcal{U} \subseteq \R^n$.
In this setting, the leader first decides on the knapsack's capacity~$x$.
Then, the uncertainties regarding the lower-level objective function
coefficients realize.
Finally, the follower solves a knapsack problem according to the
realization of his own profits, which may differ from those of the leader.
Hence, the follower decides in a wait-and-see fashion, \ie, the timing
in~\eqref{eq:wait-and-see-follower-timing} is considered.
The leader's aim is to choose the capacity of the knapsack in such
a way that her own profit of the items packed by the follower
is maximized.
Whenever the follower's choice of items is not unique, the pessimistic
approach is considered.
The deterministic variant of Problem~\eqref{eq:continuous-robust-knapsack}
can be solved in polynomial time,
which makes it a good starting point to address the
question of how uncertainties may affect the hardness of
the underlying bilevel problem.

Driven by this question, \textcite{Buchheim_Henke:2020,Buchheim_Henke:2022}
show that the complexity of Problem~\eqref{eq:continuous-robust-knapsack}
strongly depends on the considered type of the uncertainty set.
For discrete uncertainty sets as well as for interval uncertainty
under the independence assumption, \ie, for the case in which
the follower's objective function coefficients independently take values
in given intervals, Problem~\eqref{eq:continuous-robust-knapsack}
remains solvable in polynomial time.
However, the problem becomes NP-hard if the uncertainty set is
the Cartesian product of discrete sets.
In particular, this shows that replacing the uncertainty set by its
convex hull may significantly change the problem,
which is very much in contrast to the situation in single-level
robust optimization.
NP-hardness is also shown for the variants of the problem
with polytopal uncertainty sets and uncertainty sets that are defined by
a $p$-norm with~$p \in [1,\infty)$.
In particular, for all NP-hard variants of the problem, even the
evaluation of the leader's objective function is NP-hard.

As a generalization of the aforementioned works, \textcite{Buchheim_et_al:2021}
are concerned with complexity questions for robust
bilevel combinatorial problems of the form
\begin{subequations} \label{eq:robust-combinatorial-prob}
  \begin{align}
    ``\max_{x \in X} \ \min_{c \in \mathcal{U}}"
    \quad & d^\top y
            \label{eq:robust-combinatorial-obj}
    \\
    \st \quad & y \in \argmax_{y' \in \R^{n_y}} \Defset{c^\top y'
                }{
                By' \leq Ax + b}.
  \end{align}
\end{subequations}
with~$X \subseteq \set{0,1}^{n_x}$, $A \in \R^{m \times n_x}$,
$B \in \R^{m \times n_y}$, $c,\, d \in \R^{n_y}$, and~$b \in \R^m$.
Again, it is assumed that the lower-level objective function coefficients
are uncertain, that the uncertainties take values in a given uncertainty
set~$\mathcal{U} \subseteq \R^{n_y}$, and that the follower decides in a
wait-and-see fashion.
As before, the quotation marks in~\eqref{eq:robust-combinatorial-obj}
express the ambiguity in the case that the lower level does not have a
unique solution.
The deterministic variant of Problem~\eqref{eq:robust-combinatorial-prob}
is known to be NP-easy.%
\footnote{A decision problem is NP-easy if it can be polynomially
  reduced to an NP-complete decision problem \parencite{Buchheim_et_al:2021}.}
However, it is shown that interval uncertainty renders
Problem~\eqref{eq:robust-combinatorial-prob} significantly
harder than the consideration of discrete uncertainty sets.
More precisely, the robust counterpart can
be~$\Sigma^P_2$-hard\footnote{This class contains those problems that
  can be solved in nondeterministic polynomial time, provided that
  there exists an oracle that solves problems that are in NP in
  constant time.} for
interval uncertainty under the independence assumption, whereas it can be
NP-hard for uncertainty sets~$\mathcal{U}$ with~$|\mathcal{U}| = 2$ and
strongly NP-hard for general discrete uncertainty sets.
In particular, it is shown that replacing the discrete uncertainty
set by its convex hull may increase the complexity of the problem at hand,
which is in line with the results in
\textcite{Buchheim_Henke:2020,Buchheim_Henke:2022}.

\subsection{Solution Approaches for the Bilevel Knapsack Interdiction Problem with a Here-and-Now Follower}

\textcite{Beck_et_al:2022a} study discrete linear min-max problems
with uncertainties regarding the lower-level objective function
coefficients.
In contrast to the aforementioned works, which all follow the notion of strict
robustness, the authors consider a~$\Gamma$-robust approach
\parencite{Bertsimas_Sim:2003,Bertsimas_Sim:2004}.
The problem under consideration thus reads
\begin{subequations} \label{eq:Beck-et-al-prob}
  \begin{align}
    \min_x \quad & c^\top x + d ^\top y \\
    \st \quad & Ax \geq a,\, x \in X \subseteq \Z^{n_x}, \\
                 & y \in \argmax_{y' \in Y(x)} \Set{d^\top y' -
                   \max_{\defset{S \subseteq [n_y]}{|S| \leq \Gamma}}
                   \sum_{i \in S} \Delta d_iy'_i},
  \end{align}
\end{subequations}
where~$\Gamma \in [n_y] \define \set{1,\ldots,n_y}$ and~$Y(x) \subseteq \Z^{n_y}_+$
denotes the lower-level feasible set.
Here, the timing in~\eqref{eq:here-and-now-follower-timing} is considered,
\ie, both the leader and the follower decide before the uncertainty realizes.
The authors present two approaches to reformulate
Problem~\eqref{eq:Beck-et-al-prob}.
The first approach is based on an extended formulation, whereas the second one
exploits the fact that Problem~\eqref{eq:Beck-et-al-prob}
can be interpreted as a single-leader multi-follower problem with independent
followers.
Based on these reformulations, the authors propose generic branch-and-cut frameworks
to solve the problem.
Moreover, it is shown that the same techniques can also be used for the case
in which uncertainties only arise in a single packing-type constraint on the lower level.
To assess the applicability of the proposed branch-and-cut methods, the authors focus
on the~$\Gamma$-robust knapsack interdiction problem \parencite{Caprara_et_al:2016}.
In this setting, both players share a common set of items and the leader has the
ability to influence the follower's decision by prohibiting the usage of certain items
by the follower.
The authors derive problem-tailored cuts and perform a computational study on
$200$~robustified knapsack interdiction instances with up to $55$~items, \ie,
with up to $55$~variables on both the upper and the lower level.


\section{A First Glimpse at Decision Uncertainty}

Although being subject to data uncertainty, both decision makers in
the bilevel problem are assumed to take perfectly rational decisions
in the sense that they can perfectly anticipate or observe the other's
decision and that they can solve their problem to global optimality.
In decision making theory, however, it is well known that these
assumptions regarding perfect information and rationality are rarely
satisfied in a real-world context.
Luckily, bilevel optimization under uncertainty allows to relax these
assumptions in multiple ways.
Throughout this article, we assumed that the major source of
uncertainty stems from unknown or noisy input data.
However, bilevel optimization involves (at least) two decision
makers and, hence, other uncertainties in the decision making
process are also possible.
Another possible one is \emph{decision uncertainty} in which, \eg, the
leader is not sure about the reaction of the follower (for instance if
the follower does not necessarily choose an optimal solution) or in
which the follower is not sure about the observed leader's decision.
We are not going into the details here but want to give a few pointers
to the relevant literature that covers such aspects.
If the leader is uncertain about her anticipation of the
follower's optimal reaction and, thus, may want to hedge against
sub-optimal follower reactions, the resulting setup can be modeled
using so-called near-optimal robust bilevel models; see, \eg,
\textcite{Besancon_et_al:2019}.
As an extreme case of the former aspect it may be the case that
the upper-level player knows that the follower will play against
her.
This is the setting of a pessimistic bilevel optimization problem,
which is also rather naturally connected to the field of robust
optimization; see, \eg, \textcite{Wiesemann_et_al:2013}.
However, if the level of cooperation or confrontation of the
follower is not known, this leads to intermediate cases in between
of the optimistic and the pessimistic case; see, \eg,
\textcite{Aboussoror_Loridan:1995,Mallozzi_Morgan:1996}.
Moreover, in many situations it is not possible for the follower to
perfectly observe the optimal decision of the leader and the follower
thus may want to hedge against all possible leader decisions in some
uncertainty set around the observation.
Such settings are tackled in,
\eg, \textcite{Bagwell:1995,vanDamme_Hurkens:1997,Beck_Schmidt:2021}.
Finally, even if all data and the rational reaction of the follower is
known and even if the leader can, in principle, fully anticipate
the (globally) optimal reaction of the follower, it might still be the case
that limited intellectual or computational resources render it
impossible for the follower to take a globally optimal decision.
In such situations, a follower might resort to heuristic approaches
and the leader may be uncertain w.r.t.\ which heuristic is used.
For a good primer in this context, we refer to the recent paper by
\textcite{Zare_et_al:2020}.

The above list is by far not comprehensive.
A much more detailed discussion of these and other aspects can be
found in our recent survey \parencite{Beck_et_al:2022c}.
However, it is hopefully clear now how much more diverse the sources of
uncertainty can be in bilevel optimization as compared to single-level
optimization.
Hence, we expect a lot of research in this area in future years.


\printbibliography

@report{Beck_et_al:2022c,
  author = {Beck, Yasmine and Ljubić, Ivana and Schmidt, Martin},
  url    = {https://optimization-online.org/2022/06/8963/},
  date   = {2022},
  title  = {A Survey on Bilevel Optimization Under Uncertainty},
  type   = {techreport},
}

@article{Aboussoror_Loridan:1995,
  author       = {Aboussoror, Abdelmalek and Loridan, Pierre},
  url          = {https://www.semanticscholar.org/paper/Strong-weak-Stackelberg-Problems-in-Finite-Spaces-Aboussoror-Loridan/c875814905ed516f7e693d9040683c1115f31c27},
  date         = {1995},
  journaltitle = {Serdica. Mathematical Journal},
  pages        = {151--170},
  title        = {Strong-weak Stackelberg Problems in Finite Dimensional Spaces},
  volume       = {21},
}

@misc{Besancon_et_al:2019,
  author = {Besançon, Mathieu and Anjos, Miguel F. and Brotcorne, Luce},
  url    = {https://arxiv.org/pdf/1908.04040.pdf},
  date   = {2019},
  title  = {Near-optimal Robust Bilevel Optimization},
}

@misc{Buchheim_Henke:2020,
  author = {Buchheim, Christoph and Henke, Dorothee},
  url    = {https://arxiv.org/abs/1903.02810},
  date   = {2020},
  title  = {The bilevel continuous knapsack problem with uncertain follower's objective},
}

@article{Buchheim_Henke:2022,
  author =	 {Buchheim, Christoph and Henke, Dorothee},
  date =	 {2022},
  doi =		 {10.1007/s10898-021-01117-9},
  journaltitle = {Journal of Global Optimization},
  title =	 {The robust bilevel continuous knapsack problem with
                  uncertain coefficients in the follower’s objective},
  volume =	 {83},
  pages =	 {803--824},
}

@article{Buchheim_et_al:2021,
  author =	 {Buchheim, Christoph and Henke, Dorothee and
                  Hommelsheim, Felix},
  date =	 {2021},
  title =	 {On the Complexity of Robust Bilevel Optimization
                  With Uncertain Follower's Objective},
  journal =	 {Operations Research Letters},
  volume =	 {49},
  number =	 {5},
  pages =	 {703--707},
  doi =		 {10.1016/j.orl.2021.07.009},
}

@article{Kleinert_et_al:2021,
  author       = {Kleinert, Thomas and Labbé, Martine and Ljubić, Ivana and Schmidt, Martin},
  date         = {2021},
  doi          = {10.1016/j.ejco.2021.100007},
  journaltitle = {{EURO} Journal of Computational Optimization},
  pages        = {100007},
  title        = {A Survey on Mixed-Integer Programming Techniques in Bilevel Optimization},
  volume       = {9},
}

@article{Zare_et_al:2020,
  author       = {Zare, M. Hosein and Prokopyev, Oleg A. and Sauré, Denis},
  date         = {2020},
  doi          = {10.1287/deca.2019.0392},
  journaltitle = {Decision Analysis},
  number       = {1},
  title        = {On Bilevel Optimization with Inexact Follower},
  volume       = {17},
}

@article{Beck_Schmidt:2021,
  author       = {Beck, Yasmine and Schmidt, Martin},
  date         = {2021},
  doi          = {10.1016/j.orl.2021.07.010},
  issn         = {0167-6377},
  journaltitle = {Operations Research Letters},
  number       = {5},
  pages        = {752--758},
  title        = {A robust approach for modeling limited observability in bilevel optimization},
  volume       = {49},
}

@report{Beck_et_al:2022a,
  author = {Beck, Yasmine and Ljubić, Ivana and Schmidt, Martin},
  url    = {http://www.optimization-online.org/DB_FILE/2021/11/8678.pdf},
  date   = {2022},
  title  = {Exact Methods for Discrete $\Gamma$-Robust Interdiction Problems},
  type   = {techreport},
}

@book{Dempe:2002,
  author    = {Dempe, S.},
  publisher = {Springer US},
  date      = {2002},
  doi       = {10.1007/b101970},
  title     = {Foundations of Bilevel Programming},
}

@article{Wiesemann_et_al:2013,
  author       = {Wiesemann, Wolfram and Tsoukalas, Angelos and Kleniati, Polyxeni{-}Margarita and Rustem, Ber{ç}},
  date         = {2013},
  doi          = {10.1137/120864015},
  journaltitle = {SIAM Journal on Optimization},
  number       = {1},
  pages        = {353--380},
  title        = {Pessimistic Bilevel Optimization},
  volume       = {23},
}

@article{Bertsimas_Sim:2003,
  author       = {Bertsimas, Dimitris and Sim, Melvyn},
  date         = {2003},
  doi          = {10.1007/s10107-003-0396-4},
  journaltitle = {Mathematical Programming},
  pages        = {49--71},
  title        = {Robust discrete optimization and network flows},
  volume       = {98},
}

@article{Bertsimas_et_al:2011,
  author       = {Bertsimas, Dimitris and Brown, David B. and Caramanis, Constantine},
  date         = {2011},
  doi          = {10.1137/080734510},
  journaltitle = {SIAM Review},
  number       = {3},
  pages        = {464--501},
  title        = {Theory and Applications of Robust Optimization},
  volume       = {53},
}

@article{Bertsimas_Sim:2004,
  author       = {Bertsimas, Dimitris and Sim, Melvyn},
  date         = {2004},
  doi          = {10.1287/opre.1030.0065},
  journaltitle = {Operations Research},
  number       = {1},
  pages        = {35--53},
  title        = {The Price of Robustness},
  volume       = {52},
}

@article{Ben-Tal_Nemirovski:1998,
  author       = {Ben-Tal, A. and Nemirovski, A.},
  date         = {1998},
  doi          = {10.1287/moor.23.4.769},
  journaltitle = {Mathematics of Operations Research},
  number       = {4},
  pages        = {769--805},
  title        = {Robust Convex Optimization},
  volume       = {23},
}

@article{Ben-Tal_et_al:2004,
  author       = {Ben-Tal, A. and Goryashko, A. and Guslitzer, E. and Nemirovski, A.},
  date         = {2004},
  doi          = {10.1007/s10107-003-0454-y},
  journaltitle = {Mathematical Programming},
  number       = {2},
  pages        = {351--376},
  title        = {Adjustable robust solutions of uncertain linear programs},
  volume       = {99},
}

@inbook{Mallozzi_Morgan:1996,
  author    = {Mallozzi, Lina and Morgan, Jacqueline},
  editor    = {Pillo, G. Di and Giannessi, F.},
  publisher = {Springer, Boston, MA},
  booktitle = {Nonlinear Optimization and Applications},
  date      = {1996},
  doi       = {10.1007/978-1-4899-0289-4_19},
  pages     = {271--282},
  title     = {Hierarchical Systems with Weighted Reaction Set},
}

@book{Birge_Louveaux:2011,
  author    = {Birge, John R. and Louveaux, Francois},
  publisher = {Springer Science \& Business Media},
  date      = {2011},
  doi       = {10.1007/978-1-4614-0237-4},
  title     = {Introduction to Stochastic Programming},
}

@article{Hansen-et-al:1992,
  author       = {Hansen, Pierre and Jaumard, Brigitte and Savard, Gilles},
  publisher    = {SIAM},
  date         = {1992},
  doi          = {10.1137/0913069},
  journaltitle = {{SIAM} Journal on Scientific and Statistical Computing},
  number       = {5},
  pages        = {1194--1217},
  title        = {New branch-and-bound rules for linear bilevel programming},
  volume       = {13},
}

@book{Ben-Tal_et_al:2009,
  author    = {Ben-Tal, A. and El Ghaoui, L. and Nemirovski, A.},
  publisher = {Princeton University Press},
  date      = {2009},
  title     = {Robust Optimization},
}

@article{Soyster:1973,
  author       = {Soyster, A. L.},
  date         = {1973},
  doi          = {10.1287/opre.21.5.1154},
  journaltitle = {Operations Research},
  number       = {5},
  pages        = {1154--1157},
  title        = {Technical Note—Convex Programming with Set-Inclusive Constraints and Applications to Inexact Linear Programming},
  volume       = {21},
}

@article{Bagwell:1995,
  author       = {Bagwell, Kyle},
  date         = {1995},
  doi          = {10.1016/S0899-8256(05)80001-6},
  journaltitle = {Games and Economic Behavior},
  number       = {2},
  pages        = {271--280},
  title        = {Commitment and observability in games},
  volume       = {8},
}

@article{vanDamme_Hurkens:1997,
  author       = {{van Damme}, Eric and Hurkens, Sjaak},
  date         = {1997},
  doi          = {10.1006/game.1997.0524},
  journaltitle = {Games and Economic Behavior},
  number       = {1--2},
  pages        = {282--308},
  title        = {Games with Imperfectly Observable Commitment},
  volume       = {21},
}

@article{Bracken-Mcgill:1973,
  author       = {Bracken, Jerome and McGill, James T.},
  publisher    = {INFORMS},
  date         = {1973},
  doi          = {10.1287/opre.21.1.37},
  journaltitle = {Operations Research},
  number       = {1},
  pages        = {37--44},
  title        = {Mathematical programs with optimization problems in the constraints},
  volume       = {21},
}

@book{Candler-Norton:1977,
  author    = {Candler, W. and Norton, R.D.},
  publisher = {World Bank},
  url       = {https://books.google.de/books?id=TicmAQAAMAAJ},
  date      = {1977},
  series    = {Discussion Papers, Development Research Center, International Bank for Reconstruction and Development},
  title     = {Multi-level Programming},
}

@book{Stackelberg:1934,
  author    = {{von Stackelberg}, Heinrich},
  publisher = {Springer},
  date      = {1934},
  title     = {{Marktform und Gleichgewicht}},
}

@book{Stackelberg:1952,
  author    = {{von Stackelberg}, Heinrich},
  publisher = {Oxford University Press},
  date      = {1952},
  title     = {Theory of the market economy},
}

@article{jeroslow:1985,
  author       = {Jeroslow, Robert G.},
  publisher    = {Springer},
  date         = {1985},
  doi          = {10.1007/BF01586088},
  journaltitle = {Mathematical Programming},
  number       = {2},
  pages        = {146--164},
  title        = {The polynomial hierarchy and a simple model for competitive analysis},
  volume       = {32},
}

@book{Kall_Wallace:1994,
  author    = {Kall, Peter and Wallace, Stein W.},
  location  = {New York},
  publisher = {Wiley},
  date      = {1994},
  series    = {Wiley-Interscience Series in Systems and Optimization},
  title     = {Stochastic Programming},
}

@inbook{Dempe:2020,
  author    = {Dempe, Stephan},
  editor    = {Dempe, Stephan and Zemkoho, Alain},
  publisher = {Springer International Publishing},
  booktitle = {Bilevel Optimization: Advances and Next Challenges},
  date      = {2020},
  doi       = {10.1007/978-3-030-52119-6_20},
  pages     = {581--672},
  title     = {Bilevel Optimization: Theory, Algorithms, Applications and a Bibliography},
}

@inbook{Pisinger_Toth:1998,
  author    = {Pisinger, David and Toth, Paolo},
  publisher = {Springer, Boston, MA},
  booktitle = {Handbook of Combinatorial Optimization},
  date      = {1998},
  doi       = {10.1007/978-1-4613-0303-9_5},
  title     = {Knapsack Problems},
}

@incollection{DeNegre_Ralphs:2009,
  author    = {DeNegre, Scott T. and Ralphs, Ted K.},
  publisher = {Springer},
  booktitle = {Operations research and cyber-infrastructure},
  date      = {2009},
  doi       = {10.1007/978-0-387-88843-9_4},
  pages     = {65--78},
  title     = {A Branch-and-cut Algorithm for Integer Bilevel Linear Programs},
}

@article{Tang_et_al:2016,
  author       = {Tang, Y. and Richard, JP. P. and Smith, J. C.},
  date         = {2016},
  doi          = {10.1007/s10898-015-0274-7},
  journaltitle = {Journal of Global Optimization},
  pages        = {225--262},
  title        = {A class of algorithms for mixed-integer bilevel min–max optimization},
  volume       = {66},
}

@thesis{DeNegre:2011,
  author      = {DeNegre, Scott T.},
  institution = {Lehigh University},
  url         = {https://coral.ise.lehigh.edu/~ted/files/papers/ScottDeNegreDissertation11.pdf},
  date        = {2011},
  title       = {Interdiction and Discrete Bilevel Linear Programming},
  type        = {phdthesis},
}

@inproceedings{Caprara_et_al:2013,
  author    = {Caprara, Alberto and Carvalho, Margarida and Lodi, Andrea and Woeginger, Gerhard J.},
  editor    = {Goemans, Michel and Correa, José},
  publisher = {Springer, Berlin, Heidelberg},
  booktitle = {Integer Programming and Combinatorial Optimization},
  date      = {2013},
  doi       = {10.1007/978-3-642-36694-9_9},
  pages     = {98--109},
  series    = {IPCO 2013},
  title     = {A Complexity and Approximability Study of the Bilevel Knapsack Problem},
  volume    = {7801},
}

@article{Fischetti_et_al:2018b,
  author       = {Fischetti, Matteo and Monaci, Michele and Sinnl, Markus},
  date         = {2018},
  doi          = {10.1016/j.ejor.2017.11.043},
  journaltitle = {European Journal of Operational Research},
  number       = {16},
  pages        = {40--51},
  title        = {A dynamic reformulation heuristic for Generalized Interdiction Problems},
  volume       = {267},
}

@article{Fischetti_et_al:2019,
  author       = {Fischetti, Matteo and Ljubić, Ivana and Monaci, Michele and Sinnl, Markus},
  date         = {2019},
  doi          = {10.1287/ijoc.2018.0831},
  journaltitle = {INFORMS Journal on Computing},
  number       = {2},
  pages        = {390--410},
  title        = {Interdiction Games and Monotonicity, with Application to Knapsack Problems},
  volume       = {31},
}

@article{Della-Croce_Scatamacchia:2020,
  author       = {Della Croce, Federico and Scatamacchia, Rosario},
  date         = {2020},
  doi          = {10.1007/s10107-020-01482-5},
  journaltitle = {Mathematical Programming},
  pages        = {249--281},
  title        = {An exact approach for the bilevel knapsack problem with interdiction constraints and extensions},
  volume       = {183},
}

@article{Caprara_et_al:2016,
  author       = {Caprara, Alberto and Carvalho, Margarida and Lodi, Andrea and Woeginger, Gerhard J.},
  date         = {2016},
  doi          = {10.1287/ijoc.2015.0676},
  journaltitle = {INFORMS Journal on Computing},
  number       = {2},
  pages        = {319--333},
  title        = {Bilevel Knapsack with Interdiction Constraints},
  volume       = {28},
}

@article{Fortuny-Amat-Mccarl:1981,
  author       = {Fortuny-Amat, José and McCarl, Bruce},
  publisher    = {Palgrave Macmillan Journals on behalf of the Operational Research Society},
  date         = {1981},
  doi          = {10.1057/jors.1981.156},
  journaltitle = {The Journal of the Operational Research Society},
  number       = {9},
  pages        = {783--792},
  title        = {A Representation and Economic Interpretation of a Two-Level Programming Problem},
  volume       = {32},
}

\end{document}